\documentclass[a4paper, final, notitlepage, 10pt]{article}
\usepackage{latexsym,bm}
\usepackage{amsthm}
\usepackage{amsfonts}
\usepackage{amsmath}
\usepackage[all]{xy}
\usepackage[final]{graphics}
\usepackage{fancyhdr}
\usepackage{verbatim}
\usepackage{geometry}
\usepackage{xcolor}
\usepackage{enumerate}

\newtheorem{thm}{Theorem}[section]
\newtheorem{cor}[thm]{Corollary}
\newtheorem{lem}[thm]{Lemma}
\newtheorem{prop}[thm]{Proposition}

\theoremstyle{definition}

\newtheorem{rem}[thm]{Remark}

\numberwithin{equation}{section}

\newcommand{\thmref}[1]{Theorem~\ref{#1}}
\newcommand{\secref}[1]{Section~\ref{#1}}
\newcommand{\lemref}[1]{Lemma~\ref{#1}}
\newcommand{\propref}[1]{Proposition~\ref{#1}}
\newcommand{\corref}[1]{Corollary~\ref{#1}}

\renewcommand{\O}{\mathcal{O}}
\renewcommand{\L}{\mathcal{L}}

\newcommand{\PP}{\mathbb{P}}

\newcommand{\bpi}{\overline{\pi}}

\renewcommand{\bf}{\bar{f}}
\newcommand{\e}{\varepsilon}
\newcommand{\nequiv}{\stackrel{num}{\sim}}

\renewcommand{\O}{\mathcal{O}}
\renewcommand{\L}{\mathcal{L}}

\begin{document}
\title{A Characterization of Inoue Surfaces \\with $p_g=0$ and $K^2=7$}
\author{Yifan Chen \quad \quad YongJoo Shin}
\date{}
\maketitle

\begin{abstract}
Inoue constructed the first examples of smooth minimal complex surfaces of general type with $p_g=0$ and $K^2=7$.
These surfaces are finite Galois covers of the $4$-nodal cubic surface with the Galois group, the Klein group $\mathbb{Z}_2\times \mathbb{Z}_2$. For such a surface $S$, the bicanonical map of $S$ has degree $2$ and it is composed with exactly one involution in the Galois group. The divisorial part of the fixed locus of this involution consists of two irreducible components:
one is a genus $3$ curve with self-intersection number $0$ and the other is a genus $2$ curve with
self-intersection number $-1$.

Conversely, assume that $S$ is a smooth minimal complex surface of general type with $p_g=0$, $K^2=7$
and having an involution $\sigma$.
We show that, if the divisorial part of the fixed locus of $\sigma$ consists of two irreducible components $R_1$ and $R_2$,
with $g(R_1)=3, R_1^2=0, g(R_2)=2$ and $R_2^2=-1$, then the Klein group $\mathbb{Z}_2\times \mathbb{Z}_2$ acts faithfully on $S$ and $S$ is indeed an Inoue surface.
\end{abstract}

\section{Introduction}

Let $S$ be a smooth minimal complex surface of general type with $p_g=0$.
The Bogomolov-Miyaoka-Yau inequality (\cite{bogomolov, miyaoka0, yau}) yields $1\le K^2_S \le 9$.
 The bicanonical map plays an important role in the classification of these surfaces.
 For $K_S^2 \ge 2$, the bicanonical map $\varphi$ of such a surface $S$ has a surface as the image (see \cite{xiao})
 and $\varphi$ is a morphism if $K_S^2 \ge 5$ (cf.~\cite{bombieri, reider}).
  According to the results of Mendes Lopes and Pardini \cite{mendes,bicanonical1,mendespardini,mendespardini1},
    we know that  $\deg \varphi=1$ if $K_S^2=9$; $\deg \varphi=1, 2$ if $K_S^2=7,8$; $\deg \varphi=1, 2, 4$ if $K_S^2=5,6$; and $\deg \varphi=1,2,4$ if $K_S^2=3,4$ and $\varphi$ is a morphism.

It seems that when $\varphi$ has the maximal degree, the surface $S$ is characterized by a concrete example. There are many results in this direction. Mendes Lopes and Pardini \cite{mendespardini0} show that
if $K_S^2=6$ and $\deg \varphi=4$, then $S$ is a Burniat surface with $K^2=6$.
Zhang \cite{zhang} shows that if $K_S^2=5$, $\deg \varphi=4$ and $\mathrm{Im}(\varphi)$ is smooth, then
$S$ is a Burniat surface with $K^2=5$. The second named author \cite{shin} shows that if $K_S^2=4$ and $\varphi$ is a morphism of degree $4$ and
$\mathrm{Im}(\varphi)$ is smooth, then $S$ is a Burniat surface of non-nodal type with $K^2=4$.
For Burniat surfaces, see \cite{peters, burniat1}.

As mentioned above, if $K_S^2=7, 8$, the bicanonical map $\varphi$ has degree $1$ or $2$ (cf.~\cite{bicanonical1, bicanonical2}). Pardini completely classifies the surfaces with $K_S^2=8$ and $\deg \varphi=2$. These surfaces are characterized as free quotients of  products of curves (cf.~\cite[Corollary 2.3]{doubleplane}) as well as double planes (cf.~\cite[Theorem 5.1]{doubleplane}).

In this article we focus on smooth minimal surfaces $S$ of general type with $p_g=0$, $K^2=7$ and $\deg \varphi=2$.
Despite of the existence of examples, the Inoue surfaces (see \cite[Example 4.1]{bicanonical1}), and a structure theorem (see \cite[Theorem 3.2]{bicanonical2}), a complete classification of these surfaces is still out of reach.
We explain the main difficulty.  Denote by $\sigma$ the involution associated to $\varphi$ and we call it the bicanonical involution.
To understand $S$, one needs to study the quotient surface $\Sigma:=S/\sigma$ and the fixed locus $\mathrm{Fix}(\sigma)$ of $\sigma$. Mendes Lopes and Pardini \cite[Propostion~3.1]{bicanonical2} study $\Sigma$ in detail and show that $\sigma$ has $11$ isolated fixed points.
However, an explicit description of the divisorial part of $\mathrm{Fix}(\sigma)$  is still missing.
In general,  for an involution $\sigma$ on a surface of general type with $p_g=0$ and $K^2=7$, denote by $R_\sigma$ the divisorial part of $\mathrm{Fix}(\sigma)$. Lee and the second named author \cite[Table in page 3]{leeshin} describe all the possible cases for $R_\sigma$ in terms of the genus and the self-intersection number of each irreducible component of $R_\sigma$, except the case where $\sigma$ is the bicanonical involution. See also \cite{rito}.

So it is natural to first consider the case where  $R_\sigma$  has the same irreducible decomposition as the one of  an Inoue surface. Inoue surfaces \cite{inoue} are the first examples of minimal surfaces of general type with $p_g=0$ and $K^2=7$ (see also \cite{inouemfd}). In \cite[Example~4.1]{bicanonical1}, it is shown that the bicanonical map of an Inoue surface has degree $2$.
  It is shown in \cite[Section~5]{leeshin} the divisorial part of the bicanonical involution $\sigma$ has the irreducible decomposition $R_\sigma=R_1+R_2$ with $g(R_1)=3, R_1^2=0, g(R_2)=2$ and $R_2^2=-1$ (see \propref{prop:inoue}).
Conversely, we have the following theorem.

\begin{thm}\label{thm:main}
Let $S$ be a smooth minimal surface of general type with $p_g(S)=0$, $K_S^2=7$ and having an involution $\sigma$.
Assume that the divisorial part $R_\sigma$ of the fixed locus of $\sigma$ consists of two irreducible components $R_1$ and $R_2$, with $g(R_1)=3, R_1^2=0, g(R_2)=2$ and $R_2^2=-1$.
 Then the automorphism group of $S$ contains a subgroup which is isomorphic to the Klein group $\mathbb{Z}_2 \times \mathbb{Z}_2$ and $S$ is an Inoue surface.
\end{thm}

We remark that it is not hard to show $\deg \varphi=2$ and that $\sigma$ is the birational involution
(see \lemref{lem:bicanonicalinvolution}).
The main claim of the theorem is the existence of involutions on $S$ other than $\sigma$.
We briefly mention how to prove the theorem. In \secref{sec:branchfiber}, based on the results of  \cite[Propostion~3.1]{bicanonical2}, we consider  the minimal resolution $W$ of  the quotient surface $\Sigma=S/\sigma$
and study a rational fibration $\bf \colon W \rightarrow \mathbb{P}^1$.
With the assumption of \thmref{thm:main},
  we not only analyze the singular fibers of $\bf$ in detail
 (see \propref{prop:singularfibers}) but also  find a fibration of curves of genus $2$ on $W$ which induces a hyperelliptic fibration $g$ of genus $5$ on $S$ (see \propref{prop:B2fibration} and \propref{prop:pullbackH}). For these we present several lemmas about curves of genus $2$ in \secref{sec:lemmacurve}. In particular we use a topological argument to obtain \lemref{lem:galois}. The hyperelliptic fibration $g$ implies an involution $\tau$ on $S$ such that $\tau \not=\sigma$. Therefore, we obtain three commuting involutions $\sigma,\tau,\sigma\tau$ on $S$. Then $S$ is an Inoue surface by the result of \cite{commuting} (see \propref{prop:commuting}).

\paragraph{Notation and conventions}
Throughout this article, we mainly consider projective normal surfaces with at worst ordinary double points (nodes) over $\mathbb{C}$. For such a surface $X$, we use the following notation.\\
$\mathrm{Pic}(X)$: the Picard group of $X$;\\
$p_g(X)$: the geometric genus of $X$, i.e. $h^0(X,\mathcal{O}_X(K_X))$;\\
$q(X)$: the irregularity of $X$, i.e. $h^1(X,\mathcal{O}_X)$;\\
$\equiv$: a linear equivalence among divisors on $X$;\\
$\nequiv$: a numerical equivalence among divisors on $X$;\\
$\mathrm{Num}(X)$: the quotient of $\mathrm{Pic}(X)$ by $\nequiv$;\\
$\rho(X)$: the Picard number of $X$, i.e., the rank of $\mathrm{Num}(X)$;\\
$(-n)$-curve ($n \in \mathbb{N}$): a smooth irreducible rational curve with the self intersection number $-n$;
a $(-2)$-curve is also called a nodal curve.

Throughout this article, we denote by $S$ a smooth minimal surface of general type with $p_g=0$ and $K^2=7$ and by
$\varphi \colon S \dashrightarrow \mathbb{P}^7$ the bicanonical map of $S$.

\section{Known results and Inoue surfaces}\label{sec:known}
We  first recall the results of \cite{bicanonical1} and \cite{bicanonical2} on the bicanonical map of $S$.
Then we describe the bicanonical maps and bicanonical involutions of Inoue surfaces (\cite[Example~4.1]{bicanonical1}).

We list some basic properties of $S$.
Note that $S$ has irregularity $q(S)=0$ since $p_g(S)=0$ and $S$ is of general type.
The exponential cohomology sequence gives $\mathrm{Pic}(S)\cong H^2(S, \mathbb{Z})$.
Then the Noether's formula and Hodge decomposition imply $\rho(S)=3$.

An involution on $S$ is an automorphism of order $2$ on $S$.
For an involution $\sigma$ on $S$, the fixed locus $\mathrm{Fix}(\sigma)$ is a disjoint union of smooth irreducible curves and some isolated fixed points. We denote by $R_\sigma$ the divisorial part of $\mathrm{Fix}(\sigma)$ and by $k_\sigma$ the number of isolated fixed points of $\sigma$.

According to \cite{bombieri,reider}, the bicanonical map $\varphi \colon S \dashrightarrow \mathbb{P}^7$ of $S$ is a morphism. And $\deg \varphi=1$ or $2$ by \cite{bicanonical1}.
If  $\deg \varphi=2$, then $\varphi$ induces an involution on $S$.
In this case, we call this involution the bicanonical involution.
We also denote by $\Sigma$ the quotient surface of $S$ by the bicanonical involution and by $\pi \colon S \rightarrow \Sigma$ the quotient map.
Among other results, Mendes Lopes and Pardini prove the following theorem in \cite{bicanonical2}.
\begin{thm}[cf.~Proposition~3.1 and Theorem~3.3 of \cite{bicanonical2}]\label{thm:known}
Assume $\deg \varphi=2$ and let $\sigma$ be the bicanonical involution. Then
    \begin{enumerate}[\upshape (a)]
        \item $K_SR_\sigma=7$ and $k_\sigma=11$;
        \item $\Sigma$ is a rational surface with $11$ nodes and $K_\Sigma^2=-4$;
        \item there is a rational fibration $f \colon \Sigma \rightarrow \mathbb{P}^1$ such that $f\circ\pi$
            is a genus $3$ hyperelliptic fibration (see the right triangle of the diagram \eqref{diag:diagram});
        \item $K_S$ is ample.
    \end{enumerate}
\end{thm}
\proof For (a), see  \cite[Proposition~3.3 v) and Corollary~3.6 iv)]{involution} and \cite[Lemma~4.2]{manynodes}.
For (b),(c) and (d), see \cite[Proposition 3.1 and Theorem~3.3]{bicanonical2}.

\begin{rem}\label{rem:known}
 Moreover, \cite[Proposition~3.1 ii) and Theorem~3.2 iv)]{bicanonical2} describe the singular fibers of $f$ and $f\circ \pi$ explicitly. See also \cite[Remark~3.4 i)]{bicanonical2}.
 However, we do not use these results.
 For the special case we consider, we shall apply the result
of \cite{manynodes} to describe the singular fibers of $f$ even more explicitly.
 See \propref{prop:singularfibers}  and its proof.
\end{rem}

Now we briefly introduce Inoue surfaces with $p_g=0$ and $K^2=7$.
These surfaces are the first examples of surfaces of general type with $p_g=0$ and $K^2=7$.
For explicit construction see \cite{inoue} and \cite[Example~4.1]{bicanonical1}.
Here we just mention that
 an Inoue surface can be realized as a finite Galois cover of the $4$-nodal cubic surface with the Galois group $\mathbb{Z}_2\times \mathbb{Z}_2$, and that its bicanonical map  has degree $2$ and that its bicanonical involution belongs to the Galois group (cf.~see \cite[Example~4.1]{bicanonical1}).
The next proposition describes $R_\sigma$ of the bicanonical involution $\sigma$.

\begin{prop}[cf.~\cite{leeshin}]\label{prop:inoue}
Assume that $S$ is an Inoue surface. Let $\sigma$ be the bicanonical involution of $S$.
Then $R_\sigma=R_1+R_2$ with $g(R_1)=3, R_1^2=0, g(R_2)=2$  and $R_2^2=-1$.
\end{prop}
\proof Set $B_\sigma:=\pi(R_\sigma)$. If $R_\sigma$ has the irreducible decomposition $R_\sigma=R_1+\ldots+R_r$,
then $B_\sigma$ has the irreducible decomposition $B_\sigma=B_1+\ldots+B_r$, where $B_i:=\pi(R_i)$ for $1\le i \le r$.
Moreover, $R_i \cong B_i$ and $\pi^*B_i=2R_i$ and thus $R_i^2=\frac 12B_i^2$ for $1 \le i \le r$.
 In \cite[Section 5]{leeshin}, it is shown that $B_\sigma=B_1+B_2$ with $g(B_1)=3, B_1^2=0, g(B_2)=2$ and $B_2^2=-2$. The proposition follows. \qed

We now state a proposition characterizing Inoue surfaces.
\begin{prop}\label{prop:commuting}
Let $S$ be a smooth minimal surface of general type with $p_g=0$, $K^2=7$ and $\deg \varphi=2$.
Assume that there is an involution $\tau$ on $S$ other than the bicanonical involution $\sigma$. Then the subgroup $\langle \sigma, \tau \rangle$ of the automorphism group of $S$ is isomorphic to $\mathbb{Z}_2\times \mathbb{Z}_2$ and  $S$ is an Inoue surface.
\end{prop}
\proof Since $\deg \varphi=2$, the birational involutions $\sigma$ is contained in the center of the automorphism group $\mathrm{Aut}(S)$ of $S$
(see \cite[Theorem~1.2]{chenauto} for a general statement). Therefore $\langle \sigma, \tau \rangle=\{1, \sigma, \tau, \sigma\tau\}$ and $\langle \sigma, \tau \rangle\cong\mathbb{Z}_2 \times \mathbb{Z}_2$.
And since $K_SR_\sigma=7$ by \thmref{thm:known} (a), according to \cite[Theorem~2.9]{commuting}, $K_SR_\tau=5$ and $K_SR_{\sigma\tau}=5$. Then $S$ is an Inoue surface by \cite[Theorem~1.1~(a)]{commuting}.\qed

\section{The branch divisors and the singular fibers of the rational fibration}\label{sec:branchfiber}
This section and the next are devoted to the proof of \thmref{thm:main}.
So we keep the assumption in \thmref{thm:main} throughout these two sections.

\begin{lem}\label{lem:bicanonicalinvolution}
The bicanonical map $\varphi$ of $S$ has degree $2$ and $\sigma$ is the bicanonical involution.
\end{lem}
\proof By the assumption on $R_1$ and $R_2$ in \thmref{thm:main}, the adjunction formula gives $K_SR_1=4$ and $K_SR_2=3$, and thus $K_SR_\sigma=7$.
Then $k_\sigma=11$ by \cite[Lemma~4.2]{manynodes} and thus $\varphi$ is composed with $\sigma$ by
\cite[Corollary~3.6 iv)]{involution}. Since $\deg \varphi \le 2$, we conclude $\deg \varphi=2$ and that $\sigma$ is the bicanonical involution.\qed

So $\sigma$ has $11$ isolated fixed points (see \thmref{thm:known}~(a)).
Let $\e \colon V \rightarrow S$ be the blowup of these points and denote by $E_j$ ($j=0, 1, \ldots, 10$)
the corresponding exceptional divisors.
Then $\sigma$ lifts to an involution $\bar{\sigma}$ on $V$.
Denote by $W$ the quotient of $V$  by $\bar{\sigma}$, by $\bpi \colon V \rightarrow W$ the quotient map
and by $C_j$ the image of $E_j$ under $\bpi$ for $j=0, \ldots, 10$.
Then $W$ is a smooth surface
and the curves $C_j$ are  nodal curves.
The middle square of the diagram \eqref{diag:diagram} commutes
\begin{align}\label{diag:diagram}
          \xymatrix{
              & V \ar"2,1"_{\bar{h}\circ\bpi} \ar"1,3"^{\e} \ar"2,2"_{\bpi} & S \ar"2,3"^{\pi} \ar"2,4"^{f\circ\pi}& \\
           \mathbb{P}^1   & W \ar"2,1"^{\bar{h}} \ar"2,3"^{\eta}\ar@/_3mm/_{\bf=f\circ \eta}"2,4" & \Sigma \ar"2,4"^{f} & \mathbb{P}^1}
\end{align}
where $\eta$ is the minimal resolution of $\Sigma$.
In particular, $K_W\equiv \eta^*K_\Sigma$ and $W$ is a smooth rational surface with $K_W^2=-4$ by \thmref{thm:known}~(b).

Note that $\bpi$ is a smooth double cover branched along the divisor
$B_1+B_2+\sum_{j=0}^{10}C_j$, where $B_i=\bpi(\e^*R_i)$ for $i=1,2$.
Note that $B_i$ is isomorphic to $R_i$ for $i=1, 2$ and thus smooth.
Also $B_1$ and $B_2$ are disjoint,  and they are disjoint from the nodal curves $C_0, C_1, \ldots, C_{10}$.
There is an invertible sheaf $\L \in \mathrm{Pic}(W)$ such that
\begin{align}\label{eq:coveringdata}
2\L\equiv B_1+B_2+C_0+C_1+\ldots+C_{10}.
\end{align}

According to \thmref{thm:known}~(b), $\bf:=f\circ \eta \colon W \rightarrow \mathbb{P}^1$
is a rational fibration. Denote by $F$ the general fiber of $\bf$.

We will frequently refer to the following two lemmas for calculating the intersection numbers of divisors on $W$.
\begin{lem}\label{lem:B1B2}

The smooth curves $B_1$ and $B_2$ have genus $3$ and $2$ respectively.
Moreover, $B_1^2=0$, $K_WB_1=4$, $B_2^2=-2$ and  $K_WB_2=4$.

\end{lem}
\proof We have seen $B_i \cong R_i$ for $i=1, 2$ and the first assertion follows.
Note that $\bpi^*B_i=2\e^*R_i$. So $2B_i^2=4R_i^2$ and thus $B_1^2=0$, $B_2^2=-2$.
Then by the adjunction formula $K_WB_1=4$ and $K_WB_2=4$.\qed
\begin{lem}[cf.~Proposition~3.1 in \cite{involution}]\label{lem:D}
Let $D :=2K_W+B_1+B_2$. Then
    \begin{enumerate}[\upshape (a)]
        \item $\bpi^*D \equiv \e^*(2K_S)$ and $D$ is nef and big;
        \item $D^2=14, DK_W=0, DB_1=8, DB_2=6$ and $DF=4$;
        \item If $DC=0$ for an irreducible curve $C$, then $C$ is one of the $11$
              nodal curves $C_0, C_1, \ldots, C_{10}$.
    \end{enumerate}
\end{lem}
\proof
Note that $K_S=\pi^*K_\Sigma+R_1+R_2$.
We have seen $K_W=\eta^*K_\Sigma$ and $\bpi^*B_i=2\e^*R_i$ for $i=1, 2$.
Hence (a) follows from the commutativity of the square in \eqref{diag:diagram} and the fact that
$K_S$ is ample (see \thmref{thm:known}~(d)).

Then $D^2=\frac 12 (2K_S)^2=14$.
 Since $K_W^2=-4$ and $B_1B_2=0$ ($B_1$ and $B_2$ are disjoint), by \lemref{lem:B1B2} and the definition of $D$,
we have $DK_W=0$, $DB_1=8$ and $DB_2=6$.
Since  $\e(\bpi^*F)$ is a general fiber of $f\circ \pi$, by \thmref{thm:known}~(c) and the adjunction formula,
$K_S\e(\bpi^*F)=4$.
Then $DF=\frac 12 (\bpi^*D)(\bpi^*F)=\frac 12\e^*(2K_S)\bpi^*F=\frac12(2K_S)\e(\bpi^*F)=4$.

Now we prove (c). By (a), $(\bpi^*C)(\e^*(2K_S))=(\bpi^*C)(\bpi^*D)=2CD=0$.
Then $(2K_S)\e_\ast(\bpi^*C)=0$ by the projection formula and thus $\bpi^*C$ is contracted by $\e$ since
$K_S$ is ample. Recall that $\e$  contracts exactly the curves $E_0, E_1, \ldots, E_{10}$.
So $C=\bpi(E_j)=C_j$ for some $j\in \{0, \ldots, 10\}$.\qed

Observe that the nodal curves $C_0, \ldots, C_{10}$ are contained in the singular fibers of $\bf$.
Note that $W$ has Picard number $\rho(W)=14$.
It is easy to see that the $13$ divisors $K_W, F, C_0, C_1, \ldots, C_{10}$ are linearly independent in $\mathrm{Pic}(W)$.
In the proposition below,  we find a $(-1)$-curve $G$, which is contained in a fiber of $\bf$, such that
$G$ and these $13$ divisors form a basis of $\mathrm{Pic}(W)\otimes\mathbb{Q}$.
We also compute the coefficients of the divisor classes of $B_1$ and $B_2$ with respect to this basis.

\begin{prop}\label{prop:B1B2class}
After possibly renumbering the $11$ nodal curves $C_0, C_1, \ldots, C_{10}$, we have
$$B_1 \equiv -2K_W+2F\ \text{and}\ B_2 \equiv -2K_W+F+2G+C_0,$$
where $G$ is a $(-1)$-curve such that $FG=0$ and $GC_0=1$.
Moreover, $B_1G=2$ and $B_2G=1$.
\end{prop}
\proof Note that $(2K_W+B_1)^2=0$ and $K_W(2K_W+B_1)=-4$ by \lemref{lem:B1B2}.
Then we have $h^0(W, \O_W(2K_W+B_1)) \ge 3$ by the Riemann-Roch theorem.
Since $F$ is nef, $F(2K_W+B_1) \ge 0$ and thus $FB_1 \ge -2FK_W=4$.

Similarly, since $(2K_W+B_2)^2=-2$ and $K_W(2K_W+B_2)=-4$ by \lemref{lem:B1B2},
the Riemann-Roch theorem yields $h^0(W, \O_W(2K_W+B_2)) \ge 2$.
Also $F(2K_W+B_2) \ge 0$, i.e. $FB_2 \ge -2FK_W=4$.
By \lemref{lem:D}~(b), $DF=4$ and thus $F(B_1+B_2)=F(D-2K_W)=8$.
So $FB_1=FB_2=4$.

Then $F(2K_W+B_1)=0$. Since $(2K_W+B_1)^2=0$, the Zariski lemma implies $2K_W+B_1\nequiv aF$ ($a\in \mathbb{N}$) and then
$aFK_W=(2K_W+B_1)K_W$, i.e. $-2a=-4$.
So $a=2$ and $B_1 \equiv -2K_W+2F$ since a numerical equivalence is the same as a linear equivalence on any smooth rational surface.

We have seen $\dim |2K_W+B_2| \ge 1$ and $F(2K_W+B_2)=0$.
It follows that the moving part of $|2K_W+B_2|$ is composed with $|F|$, namely
$2K_W+B_2\equiv bF+\Psi$, where $b \in \mathbb{N}$ and $\Psi$ is the fixed part of $|2K_W+B_2|$.
Then
$D(bF+\Psi)=D(2K_W+B_2)$, i.e. $4b+D\Psi=6$. By \lemref{lem:D}~(a), we have $b=1$
and $2K_W+B_2\equiv F+\Psi$.

Since $\Psi\equiv B_2+2K_W-F$, by \lemref{lem:B1B2} and \lemref{lem:D}~(b), we have
 $$D\Psi=2, F\Psi=0, K_W\Psi=-2, \Psi^2=-2, B_1\Psi=4\ \text{and}\ B_2\Psi=2.$$
Since $F\Psi=0$, every irreducible component of $\Psi$ is a smooth rational curve with negative self-intersection number.
Since $K_W\Psi=-2$, there is an irreducible component $G$ of $\Psi$ such that $K_WG<0$.
It follows that $G$ is a $(-1)$-curve by the adjunction formula and then $DG>0$ by \lemref{lem:D}~(a) and (c).
We claim that $DG=1$.
Otherwise, since $DG \le D\Psi=2$ by \lemref{lem:D}~(a),
$DG=2$ and then $D(\Psi-G)=0$.
By \lemref{lem:D}~(c), $\mathrm{Supp}(\Psi-G)$ is contained in $C_0 \cup C_1 \ldots \cup C_{10}$.
Then $B_i(\Psi-G)=0$ since $B_i$ is disjoint from $C_0, C_1, \ldots, C_{10}$ for $i=1, 2$.
So $B_1G=B_1\Psi=4$ and $B_2G=B_2\Psi=2$.
Since $D\equiv 2K_W+B_1+B_2$, we have $K_WG=\frac 12(D-B_1-B_2)G=-2$.
This contradicts that $G$ is a $(-1)$-curve and thus $DG=1$.

Then $G(B_1+B_2)=G(D-2K_W)=3$. Also $GB_1=G(-2K_W+2F)=2$ and thus $GB_2=1$.
By \eqref{eq:coveringdata}, $G(C_0+C_1+\ldots+C_{10})=2\L G-3 \not=0$. In particular, $GC_j >0$ for some $j \in \{0,1,\ldots,10\}$.
After possibly renumbering the $11$ nodal curves $C_0, C_1, \ldots, C_{10}$, we may assume $j=0$ and thus $GC_0>0$.
Since $F(G+C_0)=0$, $G+C_0$ is contained in the same singular fiber of $f$ and thus $GC_0=1$.

We have shown $FG=0, DG=1, B_1G=2, B_2G=1$ and $GC_0=1$.
It remains to show $\Psi=2G+C_0$.
Note that $\Psi G=(B_2+2K_W-F)G=-1$ and $\Psi C_0=(B_2+2K_W-F)C_0=0$.
It follows that $(\Psi-2G-C_0)^2=\Psi^2-2\Psi(2G+C_0)+(2G+C_0)^2=0$.
Since $D(\Psi-2G-C_0)=0$, $\Psi=2G+C_0$ by the algebraic index theorem.\qed

\begin{cor}\label{cor:GB2}
$2G|_{B_2}\equiv K_{B_2}$.
\end{cor}
\proof
By \propref{prop:B1B2class},
$2(K_W+B_2)\equiv (-2K_W+2F)+4G+2C_0\equiv B_1+4G+2C_0$.
Since $B_1\cap B_2=\emptyset$ and $C_0\cap B_2=\emptyset$,
we have
$2K_{B_2}\equiv 2(K_W+B_2)|_{B_2}\equiv 4G|_{B_2}$.
Note that $G \not= B_2$ and $GB_2=1$.
So $G|_{B_2}$ is an effective divisor of degree $1$.
Since $g(B_2)=2$, $2G|_{B_2} \equiv K_{B_2}$ by \lemref{lem:weierstrass}.\qed

\begin{cor}\label{cor:B1fibration}
The linear system $|-2K_W+2F|$ is a base point free pencil of curves of genus~$3$.
\end{cor}
Before the proof, we remark that we can not even conclude $|-2K_W+2F| \not =\emptyset$
from the Riemman-Roch theorem: $\chi(\O_W(-2K_W+2F))=-1$.
\proof Recall that $B_1\equiv -2K_W+2F$ (see \propref{prop:B1B2class}) and $B_1$ is smooth irreducible with $B_1^2=0$ and $g(B_1)=3$ (see \lemref{lem:B1B2}).
First assume $|-K_W+F| \not = \emptyset$.
Then for $\Delta \in |-K_W+F|$, $B_1\Delta=0$ and thus $B_1$ is disjoint from $\Delta$.
Therefore $B_1$ and $2\Delta$ generate a base point free pencil of curves of genus $3$.

It suffices to prove $|-K_W+F| \not=\emptyset$.
We first show $|-K_W+F+G| \not=\emptyset$.
Note that $(-K_W+F+G)^2=1$ and $K_W(-K_W+F+G)=1$ (see \propref{prop:B1B2class}).
By Serre duality,  $h^2(W, \O_W(-K_W+F+G))=h^0(W, \O_W(2K_W-F-G))=0$.
Then $|-K_W+F+G| \not= \emptyset$ by the Riemann-Roch theorem.

Let $\Phi \in |-K_W+F+G|$.
Since $K_W+B_2\equiv -K_W+F+2G+C_0\equiv \Phi+G+C_0$ and $p_g(W)=0$,
we conclude that $\Phi \not \ge B_2$ and thus $\Phi|_{B_2}$ is an effective divisor.
We claim  $\Phi \ge G$.
Otherwise, since $G\Phi=G(-K_W+F+G)=0$, $\Phi$ and $G$ are disjoint.
Note that $$K_{B_2}=(K_W+B_2)|_{B_2}\equiv (\Phi+G+C_0)|_{B_2}=\Phi|_{B_2}+G|_{B_2}$$
 since $C_0\cap B_2=\emptyset$.
Then $\Phi|_{B_2}\equiv G|_{B_2}$ by \corref{cor:GB2}.
Since $\deg \Phi|_{B_2}=\deg G|_{B_2}=1$ and $\Phi$ is disjoint from $G$,
we have $B_2 \cong \PP^1$, a contradiction to $g(B_2)=2$.

Hence $\Phi \ge G$ and thus $\Phi-G$ is an effective divisor in $|-K_W+F|$. \qed

The next proposition describes the singular fibers of $\bf$.
\begin{prop}\label{prop:singularfibers}
Denote by $F_0$ the singular fiber of $\bf$ containing $G+C_0$.
The rational fibration $\bf$ has exactly $5$ singular fibers (possibly renumbering the $10$ nodal curves
$C_1, \ldots, C_{10}$):
          \begin{enumerate}[\upshape (a)]
            \item $F_0=2G+2C_0+2Z+C_{9}+C_{10}$, where $Z$ is a $(-2)$-curve such that $ZG=0$ and $ZC_0=ZC_{9}=ZC_{10}=1$;
            \item the other $4$ fibers are $C_{2i-1}+2\Gamma_i+C_{2i}$,
                  where $\Gamma_i$ is a $(-1)$-curve such that
                   $\Gamma_iC_{2i-1}=\Gamma_iC_{2i}=1$ for $i=1, \ldots, 4$.
          \end{enumerate}
\end{prop}
\proof
We have the following observations.
\begin{enumerate}[\upshape (i)]
    \item The $(-1)$-curve $G$ is disjoint from the $10$ nodal curves $C_1, \ldots, C_{10}$.

           Since $2G \equiv 2K_W+B_2-F-C_0$ by \propref{prop:B1B2class}, we have $GC_i=0$ for $i=1, \ldots, 10$.
    \item $F_0 \ge 2G+C_0$.

        Actually, since $(-2K_W+2F)B_2=B_1B_2=0$ by \propref{prop:B1B2class} and $B_2$ is irreducible, $B_2$ is contained in some member of $|-2K_W+2F|$ (see \corref{cor:B1fibration}).
        Then $|-2K_W+2F-B_2| \not =\emptyset$. Since $(-2K_W+2F)-B_2 \equiv F-2G-C_0$, we have $F_0 \ge 2G+C_0$.
    \item Every irreducible component of a singular fiber of $\bf$ is either a $(-1)$-curve or a nodal curve.

          It suffices to show that $-K_W$ is $\bf$-nef, which follows from \corref{cor:B1fibration}.
\end{enumerate}

Blowing down $G$ and then blowing down the image of $C_0$,
we obtain a birational morphism $\mu \colon W \rightarrow W'$,
where  $W'$ is a smooth rational surface with $K_{W'}^2=-2$ and $\rho(W')=\rho(W)-2=12$.

Denote by $p'$ the point $\mu(G+C_0)$ on $W'$.
Note that there is a fibration $f'\colon W' \rightarrow \mathbb{P}^1$ such that $\bf=f'\circ \mu$.
Set $F_0':=\mu(F_0)$. Then $F_0'$ is a fiber of $f'$ and $p' \in F_0'$.

Set $C_j':=\mu(C_j)$ for $j=1, \ldots, 10$.
Since both $G$ and $C_0$ are disjoint from the $10$ nodal curves $C_1, \ldots, C_{10}$,
we see that $C_j'$ is a nodal curve and $p' \not \in C_j'$ for $j=1, \ldots, 10$.
So $W'$ contains $10$ pairwise disjoint nodal curves.
Applying \cite[Theorem~3.3]{manynodes} and possibly renumbering the nodal curves $C_1', \ldots, C_{10}'$, we conclude that $f'$ has exactly $5$ singular fibers: $C_{2i-1}'+2\Gamma_i'+C_{2i}'$,
where $\Gamma_i'$ is a $(-1)$-curve such that $\Gamma_i'C_{2i-1}'=\Gamma_i'C_{2i}=1$ for $i=1, \ldots, 5$.
We distinguish two cases.

\noindent\paragraph{Case 1: $F_0'$ is a smooth fiber of $f'$.}
Since $p'\in F_0'$, according to (ii),
we have $F_0=\mu^*F_0'=Z+2G+C_0$,
 where $Z$ is the strict transform of $F_0'$ and $Z$ is a nodal curve with $ZG=1$ and $ZC_0=0$.
Besides $F_0$, $\bf$ has $5$ singular fibers $C_{2i-1}+2\Gamma_i+C_{2i}$, i.e., the pullback of the singular fibers of $f'$, where $\Gamma_i:=\mu^*(\Gamma_i')$ is a $(-1)$-curve for $i=1, \ldots, 5$.

Denote by $\gamma \colon B_2 \rightarrow \mathbb{P}^1$ the restriction of
$\bf \colon W \rightarrow \mathbb{P}^1$ to $B_2$
and denote by $\mathcal{R}_{\gamma}$  the ramification divisor of $\gamma$.
Then $\deg \gamma=FB_2=4$ and then $\deg \mathcal{R}_\gamma=10$ by the Riemann-Hurwitz theorem.

Since $\Gamma_i$ ($i=1,\ldots, 5$) appears with multiplicity $2$ in a singular fiber of $\bf$,
$\mathcal{R}_\gamma \ge \Gamma_i|_{B_2}$ for $i=1,\ldots, 5$.
Similarly, $\mathcal{R}_\gamma \ge G|_{B_2}$.
Note that $\deg \Gamma_i|_{B_2}=\Gamma_iB_2=\frac12(F-C_{2i-1}-C_{2i})B_2=2$ for $i=1, \ldots, 5$ and $\deg G|_{B_2}=GB_2=1$.
We have $10=\deg \mathcal{R}_\gamma \ge 2\times 5+1=11$, a contradiction.
So Case 1 does not occur.

\noindent\paragraph{Case 2: $F_0'$ is one of the singular fiber of $f'$.}
Without loss of generality, assume $F_0'=C_9'+2\Gamma_5'+C_{10}'$.
We have seen $p'\in F_0'$ and $p' \not \in C_9' \cup C_{10}'$. So $p' \in \Gamma_5'$.
By (ii) and (iii), $F_0=\mu^*F_0'=2G+2C_0+2Z+C_{9}+C_{10}$,
where $Z$ is the strict transform of $\Gamma_5'$ and $Z$ is a nodal curve such that $ZG=0$ and $ZC_0=ZC_{9}=ZC_{10}=1$.
The pullbacks of the other $4$ singular fibers of $f'$ are as described in (b) of the proposition since $p'$ does not belong to these $4$ fibers.\qed

\begin{cor}\label{prop:restrictionB2}
We have
$$F|_{B_2}\equiv 2K_{B_2},\ \ K_W|_{B_2}\equiv 2K_{B_2},\ \ B_2|_{B_2}\equiv -K_{B_2}.$$
\end{cor}
\proof First assume $F|_{B_2}\equiv 2K_{B_2}$.
Recall that $2G|_{B_2}\equiv K_{B_2}$ by \corref{cor:GB2} and $C_0 \cap B_2=\emptyset$.
Since $B_2\equiv -2K_W+F+2G+C_0$ by \propref{prop:B1B2class},
we have
$$K_{B_2}=(K_W+B_2)|_{B_2}\equiv (-K_W+F+2G+C_0)|_{B_2}\equiv -K_W|_{B_2}+3K_{B_2}.$$
Hence $K_W|_{B_2}\equiv 2K_{B_2}$.
Then $B_2|_{B_2}\equiv K_{B_2}-K_W|_{B_2}\equiv -K_{B_2}$.

We now prove $F|_{B_2}\equiv 2K_{B_2}$.
Denote by $\gamma \colon B_2 \rightarrow \mathbb{P}^1$ the restriction of
$\bf \colon W \rightarrow \mathbb{P}^1$ to $B_2$
and denote by $\mathcal{R}_{\gamma}$  the ramification divisor of $\gamma$.
Then $\deg \gamma=FB_2=4$ and then  $\deg \mathcal{R}_\gamma=10$ by the Riemann-Hurwitz theorem.

Since the multiplicity of $G+Z$ in $F_0$ is $2$ (see \propref{prop:singularfibers}),
$\mathcal{R}_\gamma \ge (G+Z)|_{B_2}$.
The same reasoning gives $\mathcal{R}_\gamma \ge \Gamma_i|_{B_2}$ for $i=1,\ldots, 4$.
Note that $B_2$ is disjoint from the nodal curves $C_0, C_1, \ldots, C_{10}$.
It follows that
$\deg (G+Z)|_{B_2}=(G+Z)B_2=\frac 12(F-2C_0-C_{10}-C_9)B_2=2$ and similarly $\deg \Gamma_i|_{B_2}=\Gamma_iB_2=2$ for $i=1, \ldots, 4$.
Since $\deg \mathcal{R}_\gamma=10$, it follows that $\mathcal{R}_\gamma = \sum_{i=1}^{4}\Gamma_i|_{B_2}+(G|_{B_2}+Z|_{B_2})$ and that $B_2$ intersects $\sum_{i=1}^{4}\Gamma_i+(G+Z)$ transversely.

We see that $\gamma$ and its five fibers $2\Gamma_i|_{B_2}$ ($i=1, 2, 3, 4$)
and $2(G+Z)|_{B_2}$ satisfy the assumption of \lemref{lem:galois}.
\lemref{lem:galois} yields
$F|_{B_2}\equiv 2K_{B_2}$.\qed

\section{The proof of the main theorem}\label{sec:pfmainthm}
We provide the complete proof of \thmref{thm:main}. We first find a genus $2$ fibration on $W$.
Recall that $G, \Gamma_i$ ($i=1,2,3,4$) and $Z$
are contained in the singular fibers of $\bf$ (see \propref{prop:singularfibers}).

\begin{prop}\label{prop:B2fibration}Let $H:=B_2+2G+C_0$.
\begin{enumerate}[\upshape (a)]
    \item The linear system $|H|$ is a base point free pencil of curves of genus $2$.
    \item For a general smooth $H \in |H|$,
          $F|_H \equiv 2K_H$ and $(\sum_{i=1}^4\Gamma_i|_H+Z|_H)\equiv 5K_H$.
\end{enumerate}
\end{prop}
\proof
By \lemref{lem:B1B2}, \propref{prop:B1B2class} and the definition of $H$, we have $HK_W=2, HB_2=0, HG=0, HF=4$ and $HC_j=0$ for $j=0, \ldots, 10$. It follows that $H^2=H(B_2+2G+C_0)=0$ and $p_a(H)=2$ by the adjunction formula.

Since $C_0\cap B_2=\emptyset$,
we have $H|_{B_2}\equiv \O_{B_2}$ by \corref{cor:GB2} and \corref{prop:restrictionB2}.
Tensoring the exact sequence
$0 \rightarrow \O_W(-B_2) \rightarrow \O_W \rightarrow \O_{B_2} \rightarrow 0$
by $\O_W(H)$,
we obtain
$$0 \rightarrow \O_W(2G+C_0) \rightarrow \O_W(H) \rightarrow \O_{B_2} \rightarrow 0.$$
It is clear that $\dim H^0(W, \O_W(2G+C_0))=1$ and
then $H^1(W, \O_W(2G+C_0))=0$ by the Riemann-Roch theorem.
The long exact sequence of cohomology groups yields $\dim |H|=1$.

To prove (a), it remains to show that $|H|$ is base point free. Since $H^2=0$, it suffices to show that $H$ is nef.
If $HC<0$ for an irreducible curve $C$, since $B_2+2G+C_0 \in |H|$, $C=B_2, G$ or $C_0$. But we have seen
$HB_2=HG=HC_0=0$, a contradiction. So $H$ is nef.

For a general $H$,
denote by $\theta$ the restriction $\bf \colon W \rightarrow \mathbb{P}^1$ to $H$
and by $\mathcal{R}_\theta$ the ramification divisor of $\theta$.
Note that $\deg \theta=HF=4$.
Also $H$ is disjoint from $G$ and $C_j$ since $HG=0$ and $HC_j=0$ for $j=0,1 \ldots, 10$.
By \propref{prop:singularfibers}, we have $HZ=2$ and $H\Gamma_i=2$ for $i=1,\ldots, 4$.
The same reasoning as the proof of \corref{prop:restrictionB2}
yields $\mathcal{R}_\theta=\sum_{i=1}^4\Gamma_i|_H+Z|_H$ and that $H$ intersects
$\sum_{i=1}^4\Gamma_i+Z$ transversely.
Then $\theta$ and its $5$ fibers $2Z|_H$ and $2\Gamma_i|_H$ ($i=1, \ldots, 4$)
satisfy the assumption of \lemref{lem:galois}.
\lemref{lem:galois} yields $F|_H\equiv 2K_H$ and $(\sum_{i=1}^4\Gamma_i+Z)|_H\equiv 5K_H$.\qed

\begin{prop}\label{prop:pullbackH}
For a general $H \in |H|$,
 $\bpi^*H$ is a hyperelliptic curve of genus $5$.
\end{prop}
\proof
Recall that $\bpi \colon V \rightarrow W$ is a double cover determined by the
branched locus $B_1+B_2+\sum_{j=0}^{10}C_j$
and the invertible sheaf $\L$, which satisfy \eqref{eq:coveringdata}.
Hence the double cover $\bpi^*H \rightarrow H$ is determined by the data
$(B_1+B_2+\sum_{j=0}^{10}C_j)|_H$ and $\L|_H$.

By \eqref{eq:coveringdata}, \propref{prop:B1B2class} and \propref{prop:singularfibers},
we have
\begin{align*}
2\L&\equiv(-2K_W+2F)+(-2K_W+F+2G+C_0)+(C_0+C_9+C_{10})+\sum_{j=1}^8C_j\\
         &=-4K_W+3F+(2G+2C_0+C_9+C_{10})+\sum_{i=1}^4(C_{2i-1}+C_{2i})\\
         &\equiv-4K_W+3F+(F-2Z)+\sum_{i=1}^4(F-2\Gamma_i)=-4K_W+8F-2(Z+\sum_{i=1}^4\Gamma_i).
\end{align*}
Since $W$ is a smooth rational surface, $\mathrm{Pic}(W)$ is torsion free and thus
$$\L \equiv -2K_W+4F-(\sum_{i=1}^4\Gamma_i+Z).$$
Since $H|_H=\O_H$ by \propref{prop:B2fibration}~(a), $K_W|_H=K_H$ by the adjunction formula.
Then $\L|_H\equiv (-2K_W+4F-(\sum_{i=1}^4\Gamma_i+Z))|_H \equiv K_H$ by \propref{prop:B2fibration}~(b).

Note that a general $H$ is disjoint from $B_2+\sum_{j=0}^{10}C_j$
since $HB_2=0$ and $HC_j=0$ for $j=0, \ldots, 10$.
Hence $(B_1+B_2+\sum_{j=0}^{10}C_j)|_H=B_1|_H$.
Since a general $H$ intersects $B_1$ transversely and $\deg B_1|_H=B_1H=B_1(B_2+2G+C_0)=4$,
we conclude that $\bpi^*H$ is irreducible and smooth.

We have shown that the double cover $\bpi^*H \rightarrow H$ and its covering data $B_1|_H$ and $\L|_H$
satisfy the assumption of  \lemref{lem:hyperelliptic} and hence $\bpi^*H$ is a hyperelliptic curve of genus $5$.\qed

Now we complete the proof of \thmref{thm:main}.
Recall the commutativity of the square in the diagram \eqref{diag:diagram}.
Denote by $\bar{h} \colon W\rightarrow \mathbb{P}^1$ the genus $2$ fibration defined by $|H|$.
Since $HC_j=0$,
 $C_j$ is contained in the fibers of $\bar{h}$ for $j=0, \ldots, 10$.
By \propref{prop:pullbackH},
$\bar{h}\circ\bpi \colon V \rightarrow \mathbb{P}^1$ is a hyperelliptic fibration of genus $5$
(see the left triangle of the diagram \eqref{diag:diagram}).
Since $\pi^*C_j=2E_j$, we see that $E_j$ is contained in the singular fibers of $\bar{h} \circ\bpi$.
Since $\bigcup_{j=0}^{10}E_j$ is the exceptional locus of $\e$,
$\bar{h} \circ\bpi$ induces a hyperelliptic fibration of genus $5$ on $S$.
Denote this fibration by  $g \colon S \rightarrow \mathbb{P}^1$.
The hyperelliptic involutions on smooth fibers of $g$
induce an involution $\tau$ on $S$ since $S$ is minimal.
Note that $\varepsilon(H)$ is a general fiber of $g$.
The quotient of $\varepsilon(H)$ by $\sigma$ is the genus $2$ curve $\eta(H)$ (see the diagram \eqref{diag:diagram}
and \propref{prop:B2fibration}~(a)), while the one by $\tau$ is a smooth rational curve.
We conclude that $\tau$ is different from $\sigma$.
Then \thmref{thm:main} follows by \propref{prop:commuting}.

\section{Lemmas on Genus Two Curves}\label{sec:lemmacurve}
Throughout this section, we denote by $C$ a smooth projective curve of genus $2$.
We omit the proofs of  \lemref{lem:bdlofdeg2} and \lemref{lem:weierstrass}.
\begin{lem}\label{lem:bdlofdeg2}
Let $L$ be an invertible sheaf on $C$ with $\deg L=2$.
Then $h^0(C, L) \ge  1$ and $h^0(C, L) \ge 2$ if and only if $L \equiv K_C$.
\end{lem}

Recall that a point $p$ on $C$ is a Weierstrass point if $2p \equiv K_C$.
\begin{lem}\label{lem:weierstrass}
\begin{enumerate}[\upshape (a)]
    \item If $p$ is a point of $C$ such that $4p \equiv 2K_C$, then $2p \equiv K_C$.
    \item The sum of the Weierstrass points of $C$ is linearly equivalent to $3K_C$.
\end{enumerate}
\end{lem}

\begin{lem}\label{lem:bidouble}
Assume that the automorphism group of $C$ contains a subgroup $G \cong \mathbb{Z}_2\times \mathbb{Z}_2$.
Denote by $q \colon C \rightarrow C/G$ the quotient map.
Then $C/G \cong \mathbb{P}^1$ and $q^*\O_{\mathbb{P}^1}(1)\equiv 2K_C$.
\end{lem}
\proof Denote by $g_1, g_2$ and $g_3$ the nontrivial elements of $G$ and
by $C_i$ the quotient of $C$ by $\langle g_i \rangle$ for $i=1, 2, 3$.
We may assume $g(C_1) \ge g(C_2) \ge g(C_3)$.
Then $g(C_1)=1$, $g(C_2)=1$, $g(C_3)=0$ and $C/G \cong \mathbb{P}^1$ by
\cite[p.~267, V.1.10.]{riemannsurfaces}.

Since $K_{C_1}\cong \O_{C_1}$, the Riemann-Hurwitz theorem shows that
$g_1$ has two fixed points $x_1$, $x_2$ and $K_C=x_1+x_2$.
Note that $g_2$ and $g_3$ permute $x_1$ and $x_2$.
So $x_1$ and $x_2$ are mapped by $q$ to a point $x \in C/G \cong \mathbb{P}^1$.
Then $2x_1+2x_2 =q^*x\equiv q^*\O_{\mathbb{P}^1}(1)$ and thus $q^*\O_{\mathbb{P}^1}(1)\equiv 2K_C$.\qed

\begin{lem}\label{lem:galois}
Let $\gamma \colon C \rightarrow \mathbb{P}^1$ be a morphism of degree $4$.
Assume that $t_1,\ldots,t_4$ and $t_5$ are $5$ distinct points of $\mathbb{P}^1$
such that
$\gamma^*(t_i)=2x_i+2y_i\ \text{and}\ x_i \not=y_i$ for $i=1,\ldots,5$.
Then $\gamma^*\O_{\mathbb{P}^1}(1)\equiv 2K_C$ and
$\sum_{i=1}^{5}(x_i+y_i) \equiv 5K_C$.
\end{lem}

\proof Set $X:=\mathbb{P}^1\setminus\{t_1, t_2, t_3, t_4, t_5\}$ and $Y:=\gamma^{-1}(X)$.
Then $Y$ with $\gamma|_Y \colon Y \rightarrow X$ is a topological covering space of $X$.
Fix $t_0 \in X$ and take a simple loop $l_i$ in $X$ based at $t_0$ and going around $t_i$ for $i=1, \ldots, 5$.
Denote by $[l_i]$ the class of $l_i$ in $\pi_1(X, t_0)$.
Then $\pi_1(X, t_0)$ is generated by $[l_1],\ldots, [l_5]$ with the relation $[l_1][l_2][l_3][l_4][l_5]=1$.

Identify the permutation group of $\gamma^{-1}(t_0)$
with the symmetric group $\mathrm{S}_4$ with $4$ letters.
The group $\pi_1(X,t_0)$ acts on $\gamma^{-1}(t_0)$ (from the right) and
corresponds to an anti-group homomorphism
 $\alpha \colon \pi_1(X, t_0) \rightarrow \mathrm{S}_4$.
 Also the group $\mathrm{D}$ of Deck transformations of the covering $Y \rightarrow X$ acts on $\gamma^{-1}(t_0)$ (from the left)
 and corresponds to a group homomorphism $\beta \colon \mathrm{D} \rightarrow \mathrm{S}_4$.
 It is well known that $\beta$ is injective and $\mathrm{Im}(\beta)$ is the centralizer of $\mathrm{Im}(\alpha)$ in
 $\mathrm{S}_4$.

For $i=1,\ldots, 5$, since $\gamma^*(t_i)=2x_i+2y_i$ and $x_i \not=y_i$,
 we conclude that $\alpha([l_i]) \not =1$ and $\alpha([l_i]) \in \mathrm{V}_4$,
  where $\mathrm{V}_4:=\{1, (12)(34), (13)(24), (14)(23)\}.$
 Since any two nontrivial elements of $\mathrm{V}_4$ generate $\mathrm{V}_4$ and
 $\alpha([l_5]) \dots \alpha([l_1])=\alpha([l_1] \dots [l_5])=1$,
 we have $\mathrm{Im}(\alpha)=\mathrm{V}_4$.
Since the centralizer of $\mathrm{V}_4$ in $\mathrm{S}_4$ is itself,
$\mathrm{D} \cong \mathrm{Im}(\beta) =\mathrm{V}_4 \cong \mathbb{Z}_2\times \mathbb{Z}_2$.

Note that $\mathrm{D}$ is indeed isomorphic to the Galois group of $\gamma$.
So $\gamma^*\O_{\mathbb{P}^1}(1) \equiv 2K_C$ by \lemref{lem:bidouble}.
The Riemann-Hurwitz theorem yields
$K_C =\gamma^*K_{\mathbb{P}^1}+\sum_{i=1}^5(x_i+y_i)$.
Since $K_{\mathbb{P}^1}\equiv \O_{\mathbb{P}^1}(-2)$,
we have $\sum_{i=1}^5(x_i+y_i)\equiv5K_C$.\qed

\begin{lem}\label{lem:hyperelliptic}
Let $\pi \colon E \rightarrow C$ be a double cover from a smooth projective curve $E$ onto $C$.
 Assume that the branched locus of $\pi$ consists of $4$ points $x_1, x_2, x_3, x_4$ such
 that $x_1+x_2+x_3+x_4 \equiv 2K_C$ and $\pi_\ast\O_E \cong \O_C\oplus \O_C(-K_C)$.
Then $E$ is a hyperelliptic curve of genus $5$.
\end{lem}
\proof The Riemann-Hurwitz theorem yields $g(E)=5$.
It suffices to show that the canonical image of $E$ is a rational curve in $\mathbb{P}^4$.

Note that $|2K_C|$ is composed with the hyperelliptic pencil $|K_C|$.
We may assume $x_1+x_2, x_3+x_4 \in |K_C|$.
Choose $s_1, s_2 \in H^0(C, \O_C(K_C))$ such that
$(s_1)_0=x_1+x_2, (s_2)_0=x_3+x_4.$
Then $s_1, s_2$ is a basis of $H^0(C, \O_C(K_C))$ and thus
$s_1^2, s_1s_2, s_2^2$ is a basis of $H^0(C, \O_C(2K_C))$.

Denote by $s_j'$ the pullback of $s_j$ by $\pi$ for $j=1, 2$,
and by $y_i$ the inverse image of $x_i$ for $i=1, 2, 3, 4$.
Let $s' \in H^0(E, \O_E(y_1+y_2+y_3+y_4))$ such that
$(s')_0=y_1+y_2+y_3+y_4$.
Since
$2y_1+2y_2+2y_3+2y_4=\pi^*(x_1+x_2+x_3+x_4)$,
we may choose $s'$ such that $s'^2=s_1's_2'$.

The assumption $\pi_\ast\O_E\cong \O_C\oplus \O_C(-K_C)$ implies
\begin{align*}
H^0(E, \O_E(K_E))\cong s'\pi^*H^0(C, \O_C(K_C)) \oplus \pi^*H^0(C, \O_C(2K_C))
\end{align*}
(cf.~\cite[Proposition~4.1]{pardini}).
It follows that
$s's_1', s's_2'$ and $s_1'^2, s_1's_2', s_2'^2$ together form  a  basis of $H^0(E, \O_E(K_E))$.
Since $s'^2=s_1's_2'$, the image of $E$ under the map $E \rightarrow \mathbb{P}^4$ defined by the basis $s_1'^2, s's_1', s_1's_2', s's_2', s_2'^2$
satisfies the equations
$z_iz_j-z_lz_k=0$
for $1 \le i, j, k, l \le 5$ and $i+j=l+k$,
where $[z_1, z_2, z_3, z_4, z_5]$ is the homogeneous coordinates of $\mathbb{P}^4$.
Hence the canonical image of $E$ is a rational normal curve. \qed

\paragraph{Acknowledgement.}
The first named author is greatly indebted to Yi~Gu for many discussions.
 The first named author would like to thank Meng Chen for the invitation to Fudan University, Wenfei Liu for
 the invitation to Xiamen University and for their hospitality. The second named author would like to thank Seonja Kim for useful comments of curves.

The first named author was supported by the National Natural Science Foundation of China (Grant No.:~11501019). The second named author was supported by Basic Science Research Program through the National Research Foundation of Korea (NRF) funded by the Ministry of Education (No. 2017R1D1A1B03028273).

\noindent\textbf{Authors' Addresses:}\\\smallskip

\noindent Yifan~Chen,\\
School of Mathematics and Systems Science, Beijing University of Aeronautics and Astronautics, \\
Xueyuan Road No.~37, Beijing 100191, P.~R.~China\\
Email:~chenyifan1984@gmail.com\\\smallskip

\noindent YongJoo Shin,\\
Department of Mathematical Sciences, KAIST, 291 Daehak-ro, Yuseong-gu, Daejeon 34141, Republic of Korea\\
Email: haushin@kaist.ac.kr

\end{document}